\RequirePackage{iftex}
\ifXeTeX
  \RequirePackage{scrlfile}
  \PreventPackageFromLoading{breakurl}
\fi
\documentclass[hidelinks,onefignum,onetabnum,nohypdvips]{siamart251216}
\usepackage{amssymb,mathtools}
\usepackage{array,booktabs,tabularx,placeins,microtype}
\usepackage{needspace}
\headers{A counterexample to a convergence bound}{J. Xu and J. Hu}
\hypersetup{
  pdftitle={A counterexample to a three-step convergence bound for BFGS with exact line search},
  pdfauthor={Jiawei Xu and Jiang Hu},
  pdfkeywords={BFGS, exact line search, convergence rates, counterexample}
}
\newcommand{\R}{\mathbb R}
\newcommand{\norm}[1]{\left\|#1\right\|}
\newcommand{\U}{\mathcal U}
\newcommand{\Sstep}{\mathcal S}
\DeclareMathOperator{\diag}{diag}
\DeclareMathOperator{\Span}{span}

\title{A counterexample to a three-step convergence bound
for BFGS with exact line search}
\author{
Jiawei Xu\thanks{Qiuzhen College, Tsinghua University
(\email{xu-jw23@mails.tsinghua.edu.cn}).}
\and
Jiang Hu\thanks{Yau Mathematical Sciences Center, Tsinghua University
(\email{jianghu@tsinghua.edu.cn}).}
}

\begin{document}
\maketitle

\begin{abstract}
We disprove a three-step error estimate for the
Broyden--Fletcher--Goldfarb--Shanno {(BFGS)} method
with exact line searches.
We construct a globally strongly convex function
$F\in C^\infty(\mathbb R^3)$ and a single nonterminating trajectory
$x_k\to0$ for which
$\limsup_{k\to\infty}\|x_{k+3}\|/(\|x_{k+2}\|\|x_k\|)=+\infty$.
The Hessian of $F$ lies between $\tfrac12 I$ and $\tfrac32 I$
everywhere, and $F(x)-\tfrac12\|x\|^2$ is compactly supported
with all derivatives vanishing at $0$. The proof first constructs
point-gradient pairs satisfying the exact BFGS step identities.
A five-step cycle produces a large three-step quotient and restores
an invariant plane for the actual updated matrix, allowing the
cycle to be repeated with arbitrarily small gradient errors.
Disjoint smooth corrections to a quadratic then realize all the
pairs by one objective. Strong convexity identifies the prescribed
points with the unique exact line-search minimizers.
Thus smoothness and strong convexity do not imply the proposed
product estimate, even for an objective that agrees with a
quadratic to every order at its minimizer.
\end{abstract}



\section{Introduction}\label{sec:introduction}

The Broyden--Fletcher--Goldfarb--Shanno (BFGS) method is one of the
most widely used methods for smooth unconstrained optimization,
owing to its practical efficiency and its ability to approximate
curvature using only gradient information
\cite[Chapter~6]{NocedalWright2006}.
Its finite termination on strictly convex quadratics under exact
line searches motivates the search for quantitative error estimates
over several successive iterations on smooth nonlinear objectives
\cite[Section~2]{Powell1984}.
Let $(x_k)\subset\R^3$ be a non-terminating BFGS trajectory generated
with exact line searches and converging
to a minimizer $x^*$ of a smooth function, with
$\nabla^2F(x^*)\succ0$. The estimate of interest is
\begin{equation*}\tag{P}\label{eq:product}
\norm{x_{k+3}-x^*}
=O\!\left(\norm{x_{k+2}-x^*}\,\norm{x_k-x^*}\right).
\end{equation*}
The implied constant may depend on the objective and initial data,
but must be independent of $k$ along this one trajectory.
Product estimates of this type were discussed by Powell
\cite[equation~(1.12) and Section~7]{Powell1984}.
The question is whether smoothness and positive
definite curvature are sufficient for \eqref{eq:product}, or whether
additional conditions on the successive search directions are needed.

Variable metric methods began with the
Davidon--Fletcher--Powell (DFP) method
\cite{Davidon1991,FletcherPowell1963}, which updates an inverse
Hessian approximation using gradient information.
BFGS arose from the subsequent updates of
Broyden, Fletcher, Goldfarb and Shanno
\cite{Broyden1970,Fletcher1970,Goldfarb1970,Shanno1970}.
Dennis and Mor\'e \cite{DennisMore1977} give a systematic account
of the derivation and convergence theory of quasi-Newton methods.

Powell \cite{Powell1971} analyzed variable metric methods
with exact line searches, and Dixon \cite{Dixon1972} established
conditions under which different members of the Broyden family
generate identical iterates.
Broyden, Dennis and Mor\'e \cite{BroydenDennisMore1973} proved local
superlinear convergence without line searches for sufficiently
accurate initial points and matrix approximations.
Dennis and Mor\'e \cite{DennisMore1974} characterized superlinear
convergence by the accuracy of the derivative approximation along
the steps. Under Wolfe line-search conditions
\cite{Wolfe1969,Wolfe1971} and appropriate regularity assumptions,
Byrd, Nocedal and Yuan \cite{ByrdNocedalYuan1987} established
global and superlinear convergence for the restricted Broyden
class excluding DFP on convex problems. Byrd and
Nocedal \cite{ByrdNocedal1989} developed further tools for analyzing
the matrix update on uniformly convex problems.

On a strictly convex quadratic in $n$ dimensions, exact-line-search BFGS terminates
in at most $n$ steps \cite[Section~2]{Powell1984}. The step rule
matters even on quadratics: Powell \cite{Powell1986} analyzed poor
transient behavior of BFGS and DFP with fixed unit steps.
For exact line searches near a nonquadratic minimizer, finite
termination on quadratics suggests bounds between
errors separated by $n$ iterations. Ritter \cite{Ritter1980}
obtained strengthened rate estimates under additional conditions.
Powell \cite[Sections~1 and~7]{Powell1984} discussed the role of
conditions on successive directions and the extent to which such
estimates could hold without them, treating the two-dimensional
case separately. Our result concerns precisely \eqref{eq:product}.

Rodomanov and Nesterov \cite{RodomanovNesterovGreedy2021} obtained
explicit superlinear bounds for greedy quasi-Newton methods,
which select update directions separately from the steps.
For classical methods with unit steps, Rodomanov and Nesterov
\cite{RodomanovNesterov2022,RodomanovNesterovNew2021} and Jin and
Mokhtari \cite{JinMokhtari2023} give explicit local superlinear
bounds under regularity and initialization assumptions.
Jin, Jiang and Mokhtari analyze BFGS with exact line searches
\cite{JinJiangMokhtari2026} and with Armijo--Wolfe line searches
\cite{JinJiangMokhtari2024}. They establish global linear bounds
for strongly convex objectives with Lipschitz continuous gradients
and superlinear bounds under additional Hessian regularity,
allowing arbitrary initial points and positive definite initial
matrices. These bounds estimate errors after a prescribed number
of iterations. They do not directly control the quotient in
\eqref{eq:product}, whose denominator contains two actual errors
from the same trajectory.

Adjacent-iterate convergence order in the quotient
($q$) sense is a distinct issue. Yuan
\cite{Yuan1984} constructed twice continuously differentiable,
uniformly convex examples with $q$-order $1+1/N$ for arbitrarily
large $N$, showing that the infimum of the possible orders is one
(see also Dai \cite{Dai2022}).
Liu, Li and Wen \cite{LiuLiWen2026} attain this boundary with a
single smooth globally strongly convex objective: exact-line-search
BFGS converges $q$-superlinearly, but no fixed $p>0$ yields an
eventual bound $\norm{x_{k+1}-x^*}\le C\norm{x_k-x^*}^{1+p}$.
Estimate \eqref{eq:product} instead couples the
errors at $k$, $k+2$ and $k+3$ along one infinite run.

Counterexamples to convergence provide another relevant comparison.
Dai \cite{Dai2002} constructed a six-point cycling example showing
that Wolfe line searches do not ensure convergence for nonconvex
objectives. Mascarenhas \cite{Mascarenhas2004} showed that failure
can also occur with exact line searches. Dai's later
four-dimensional polynomial example
\cite{Dai2013} has unit steps that uniquely minimize every line-search
function, while the BFGS run fails to converge. Modified updates
can recover global convergence under suitable assumptions in
nonconvex minimization; Li and Fukushima \cite{LiFukushima2001}
provide one such construction. For uniformly convex
objectives, Yuan \cite{Yuan1995} obtained DFP convergence results under
additional conditions on the generated sequence. His survey
\cite{Yuan1999} recorded convergence under Wolfe line searches without
such additional conditions as an open question. More recently, Liu
et al.\ \cite{LiuWangWenZhangYuan2026} constructed a $C^2$ uniformly
convex objective in $\R^2$ and a classical DFP run whose accepted
steps satisfy the standard strong Wolfe conditions, while the
gradient norms converge to a positive constant.
For \eqref{eq:product},
the issue is how rapidly a trajectory approaches a nondegenerate
minimizer once convergence is known.

Our contributions are as follows.
\begin{enumerate}
\item We disprove \eqref{eq:product} for a globally strongly convex
$F\in C^\infty(\R^3)$ with
$\tfrac12 I\preceq\nabla^2F\preceq\tfrac32 I$. Along one infinite
BFGS trajectory converging to $0$, the product quotients have
infinite limsup, even though $F(x)-\tfrac12\norm{x}^2$ is compactly
supported and all its derivatives vanish at $0$.
\item We construct the point-gradient pairs by repeating a five-step
cycle. Three steps produce a large quotient, and a further gradient
change restores the invariant plane needed for the next cycle.
This gives unbounded quotients within one sequence while allowing
arbitrarily small gradient errors.
\item We construct one smooth function with all the prescribed
gradients by adding disjoint local corrections to the quadratic $\tfrac12\norm{x}^2$.
The construction preserves the Hessian bounds and smoothness at
$0$. Strong convexity then makes the prescribed sequence exactly
the BFGS trajectory with exact line search.
\end{enumerate}

Section~\ref{sec:result} gives the BFGS iteration, the main theorem,
and the two key lemmas. Section~\ref{sec:tools} constructs the local
pairs; Section~\ref{sec:cycle} combines Table~\ref{tab:cycle} with
the quotient and continuation lemmas.
Section~\ref{sec:objective} concatenates the cycles, constructs
the objective, and identifies its actual BFGS trajectory.

\Needspace{14\baselineskip}
\section{BFGS and the main result}\label{sec:result}

Given an objective $F\in C^3(\R^3)$, an initial point $x_1\in\R^3$,
and a symmetric positive definite matrix $B_1$, BFGS with exact
line search takes the form
\begin{equation}\label{eq:bfgs}
\begin{gathered}
g_k:=\nabla F(x_k),\qquad d_k:=-B_k^{-1}g_k,\\
\alpha_k\in\operatorname*{argmin}_{\alpha\in\R}F(x_k+\alpha d_k),
\qquad x_{k+1}:=x_k+\alpha_kd_k,\\
s_k:=x_{k+1}-x_k,\qquad y_k:=g_{k+1}-g_k,\\
B_{k+1}:=B_k-\frac{B_ks_ks_k^TB_k}{s_k^TB_ks_k}
              +\frac{y_ky_k^T}{s_k^Ty_k}.
\end{gathered}
\end{equation}
Here the line search is over the entire affine line. In the
counterexample below, every line minimizer is unique and has a
positive parameter, and every curvature product $s_k^Ty_k$ is
positive.

All vector norms are Euclidean, and matrix norms are induced
operator norms. We write $I=I_3$ and use $e_1,e_2,e_3$ for the
standard basis. For symmetric matrices, $A\preceq C$ means that
$C-A$ is positive semidefinite. To construct $F$, we will add small
smooth corrections to the reference quadratic
\[
Q(x):=\tfrac12\norm{x}^2,\qquad \nabla Q(x)=x,\qquad \nabla^2Q(x)=I.
\]
A smooth function is \emph{flat at $0$} if all of its derivatives,
including its value, vanish there. Our main result is the following.

\begin{theorem}\label{thm:main}
There exist $F\in C^\infty(\R^3)$, $x_1\in\R^3\setminus\{0\}$, and
$B_1=B_1^T\succ0$ with the following properties.
The function $F-Q$ is compactly supported and flat at $0$, and
\[
\tfrac12I\preceq\nabla^2F(x)\preceq\tfrac32I
\qquad(x\in\R^3).
\]
$0$ is the unique minimizer of $F$, with $\nabla^2F(0)=I$, and
$\nabla F(x)\ne0$ for $x\ne0$.
The BFGS iteration \eqref{eq:bfgs} is well defined for every $k\ge1$:
each line search has a unique minimizing parameter $\alpha_k>0$,
each matrix $B_k$ is positive definite, and
$x_k\ne0$, $g_k\ne0$ and $s_k^Ty_k>0$.
Its iterates converge to $0$ and, for every integer $m\ge1$,
\begin{equation}\label{eq:main-quotient}
\frac{\norm{x_{5m}}}{\norm{x_{5m-1}}\norm{x_{5m-3}}}\ge2m.
\end{equation}
Consequently, the product quotient in \eqref{eq:product} has infinite
limsup and \eqref{eq:product} is false.
\end{theorem}

\subsection{Two lemmas for the proof}\label{sec:realization}

The proof has two parts. We first construct a sequence of
point-gradient pairs $(x_j,g_j)$ whose points violate the product
estimate and whose data satisfy the BFGS step identities. We then
realize these pairs by one smooth strongly convex objective and
show that its actual BFGS run follows exactly the prescribed
sequence. Until this realization, $g_j$ is a prescribed vector;
it is not yet the gradient of a specified function.

For use in the discrete construction, denote the BFGS update by
\begin{equation}\label{eq:update}
\U(B;s,y):=B-\frac{Bss^TB}{s^TBs}+\frac{yy^T}{s^Ty},
\qquad B=B^T\succ0,\quad s^Ty>0.
\end{equation}

\begin{lemma}[Construction of the points and gradients]\label{lem:trajectory}
For every $c>0$, there exist sequences $x_j,g_j\in\R^3$,
symmetric positive definite matrices $B_j$, and numbers $\beta_j>0$,
with $g_j\ne0$ and $0<\norm{x_1}<1/2$, such that the following
properties hold. Writing $s_j:=x_{j+1}-x_j$ and
$y_j:=g_{j+1}-g_j$, we have, for every $j\ge1$,
\begin{equation}\label{eq:algebraic-run}
s_j=-\beta_jB_j^{-1}g_j,\qquad
g_{j+1}^Ts_j=0,\qquad s_j^Ty_j>0,\qquad
B_{j+1}=\U(B_j;s_j,y_j),
\end{equation}
and
\begin{equation}\label{eq:flat-data}
0<\norm{x_{j+1}}<\norm{x_j}/8,\qquad
\norm{g_j-x_j}\le c\,2^{-j}\norm{x_j}^{j+3}.
\end{equation}
Moreover, for every integer $m\ge1$,
\begin{equation}\label{eq:bad-subsequence}
\frac{\norm{x_{5m}}}{\norm{x_{5m-1}}\norm{x_{5m-3}}}\ge2m.
\end{equation}
In particular, $x_j\to0$, and the three-step product estimate
\eqref{eq:product} fails along this single sequence.
\end{lemma}

Lemma~\ref{lem:trajectory} is proved in Section~\ref{sec:trajectory}
using the five-step construction in Section~\ref{sec:cycle}. The first bound in \eqref{eq:flat-data}
separates the points and forces convergence to $0$. The second
makes the gradient errors smaller than every fixed power of the
point radius as the index increases. These bounds allow the
following realization result.

\begin{lemma}[A function with the prescribed BFGS iterates]\label{lem:realization}
There exists $c_0>0$ with the following property.
Let $(x_j,g_j,B_j,\beta_j)_{j\ge1}$ satisfy the conclusions of
Lemma~\ref{lem:trajectory} with $c=c_0$, including
\eqref{eq:algebraic-run}--\eqref{eq:bad-subsequence}.
Then there exists $F\in C^\infty(\R^3)$ such that $F-Q$ is
compactly supported and flat at $0$, and
\begin{equation}\label{eq:realization}
F(x_j)=Q(x_j),\qquad \nabla F(x_j)=g_j,\qquad
\tfrac12I\preceq\nabla^2F(x)\preceq\tfrac32I
\end{equation}
for all $j\ge1$ and $x\in\R^3$.
$0$ is the unique minimizer of $F$, with $\nabla^2F(0)=I$, and
$\nabla F(x)\ne0$ for $x\ne0$.
The BFGS iteration \eqref{eq:bfgs} applied to $F$ from $(x_1,B_1)$
produces exactly $(x_j,g_j,B_j)$, with unique line-search parameters
$\alpha_j=\beta_j>0$. Its iterates therefore satisfy
\eqref{eq:bad-subsequence}.
\end{lemma}

The proof of Lemma~\ref{lem:realization} is given in
Section~\ref{sec:smooth-objective},
with derivative estimates in Appendix~\ref{app:interpolation}. Disjoint smooth corrections to
$Q$ interpolate the gradients while preserving the Hessian bounds.
Global strong convexity and the orthogonality in
\eqref{eq:algebraic-run} then identify every prescribed point
as the unique exact line-search minimizer. In particular, $c_0$
is fixed before the discrete trajectory is constructed.

\begin{proof}[Proof of Theorem~\ref{thm:main}]
Choose the constant $c_0$ from Lemma~\ref{lem:realization} and
apply Lemma~\ref{lem:trajectory} with $c=c_0$.
Lemma~\ref{lem:realization} supplies a single objective $F$ with
the required smoothness, support, flatness, Hessian bounds, and
minimizer, and identifies all the constructed data with its
actual BFGS run. The step identities in \eqref{eq:algebraic-run}
give positive curvature products, and \eqref{eq:flat-data} gives
nonzero iterates converging to $0$. Every $B_j$ is positive
definite and every $g_j$ is nonzero by Lemma~\ref{lem:trajectory}.

Finally, \eqref{eq:bad-subsequence} is exactly
\eqref{eq:main-quotient}. Taking $k=5m-3$ shows that the quotient
in \eqref{eq:product} is at least $2m$ along this subsequence.
All its denominators are positive, so its limsup is infinite.
\end{proof}

Most prescribed pairs in Lemma~\ref{lem:trajectory} satisfy $g_j=x_j$;
only two gradients in each five-step cycle are changed. The next
sections construct those changes with the arbitrarily small
tolerances required in \eqref{eq:flat-data}.

\section{Construction of the points and gradients}\label{sec:tools}

The construction has two parts: a small gradient change creates a
large three-step quotient, and a second change prepares the next
cycle. All other prescribed gradients equal the corresponding
points. We first record the conditions for one step, then give
the two lemmas that make this construction repeatable.

A triple $(z,g,B)$ consists of a point, its prescribed gradient,
and its current positive definite matrix. For a step $z\to z_+$,
the \emph{endpoint} is the next point $z_+$ and the displacement
is $z_+-z$. We use $z_i$ in local coordinates and $x_j$ for the
full sequence in fixed coordinates.

\subsection{Properties of one BFGS step}

When constructing the next gradient, we must ensure that it is
orthogonal to the step. The following lemma states this condition
together with the required search direction and step length.

\begin{lemma}[Properties of one BFGS step]\label{lem:transition}
Let $(z,g)$ and $(z_+,g_+)$ be two point-gradient pairs, with
$g\ne0$, and let $B=B^T\succ0$. Set $s:=z_+-z$ and
$y:=g_+-g$. Suppose
\begin{equation}\label{eq:step-conditions}
z_+-z=-\alpha B^{-1}g,\qquad \alpha>0,\qquad
g_+^T(z_+-z)=0.
\end{equation}
Then
\[
s^Ty=\alpha g^TB^{-1}g>0.
\]
The updated matrix $C:=\U(B;s,y)$ is positive definite and
satisfies the secant equation $Cs=y$.
If a differentiable strongly convex function $F$ satisfies
$\nabla F(z)=g$ and $\nabla F(z_+)=g_+$, then $z_+$ is the
unique minimizer of $F$ on the search line $z-tB^{-1}g$, $t\in\R$.
The step conditions and matrix update are preserved
under orthogonal changes of coordinates.
\end{lemma}

\begin{proof}
The step conditions give $s\ne0$ and
$s^Ty=s^T(g_+-g)=-s^Tg=\alpha g^TB^{-1}g>0$.
Positive definiteness and the secant equation follow from the
standard BFGS properties \cite[Section~6.1]{NocedalWright2006}.
For such a function $F$, set $\phi(t):=F(z-tB^{-1}g)$. The orthogonality condition gives
$\phi'(\alpha)=-g_+^TB^{-1}g=0$. Thus $\alpha$ is the unique
minimizer of the strongly convex function $\phi$.
For an orthogonal matrix $O$,
\[
(O^TBO)^{-1}O^Tg=O^TB^{-1}g,\qquad
\U(O^TBO;O^Ts,O^Ty)=O^T\U(B;s,y)O.
\]
Inner products and the line parameter are preserved, proving the
last assertion.
\end{proof}

\subsection{Choosing the next point and gradient}

Given $(z,g,B)$, we use the direction $-B^{-1}g$, minimize $Q$
on that line, and prescribe the next gradient as
$g_+=z_+=\nabla Q(z_+)$. For data in
\[
\mathcal A:=\{(z,g,B):B=B^T\succ0,\ z^TB^{-1}g>0\},
\]
define $\Sstep(z,g,B):=(z_+,C)$ by
\begin{equation}\label{eq:transition}
w:=B^{-1}g,\qquad
\lambda:=\frac{z^Tw}{w^Tw},\qquad
z_+:=z-\lambda w,\qquad
C:=\U(B;z_+-z,z_+-g).
\end{equation}
The output triple is $(z_+,z_+,C)$. Minimizing $Q$ gives
\eqref{eq:step-conditions} and
\begin{equation}\label{eq:compatibility}
z_+^T(z_+-z)=0,\qquad
(z_+-z)^T(z_+-g)=\lambda g^TB^{-1}g>0.
\end{equation}
When $g=z$, this is one BFGS step for $Q$; for general $g$, the
search direction uses the prescribed gradient instead. Write
\begin{equation}\label{eq:quadratic-endpoint}
\lambda_Q(z,B):=\frac{z^TB^{-1}z}{\norm{B^{-1}z}^2},\qquad
Z_Q(z,B):=z-\lambda_Q(z,B)B^{-1}z
\quad(z\ne0,\ B\succ0).
\end{equation}
so that $Z_Q(z,B)$ is the minimizer of $Q$ on the search line;
it equals $0$ exactly when $z$ is an eigenvector of $B$.

The first step of a cycle starts from $z_{-1}$, with prescribed
gradient $z_{-1}$ and matrix $D$. We use the following form:
\begin{equation}\label{eq:canonical}
\begin{gathered}
z_{-1}=(L,R,0)^T,\qquad z_0=(0,R,0)^T,\\
D=\frac{\theta}{L^2}z_{-1}z_{-1}^T+\diag(0,a,b),
\qquad a,b>0.
\end{gathered}
\end{equation}
Here $L$ is the step length, $R=\norm{z_0}$, and $\theta$ is
the line parameter. The coefficients $a,b$ remain after the
rank-one term in $D$ is removed by the BFGS update.
The next lemma identifies the point reached by this step and
the updated matrix after choosing its gradient. The gradient
must be perpendicular to the displacement, so it has the form
$g=(0,R+T,H)^T$.

\begin{lemma}[The first step and its updated matrix]\label{lem:endpoint-update}
For the data in \eqref{eq:canonical}, with $L,R,\theta,a,b>0$,
the exact line search for $Q$ from $(z_{-1},z_{-1},D)$ reaches
$z_0$ with parameter $\theta$. Prescribe $g=(0,R+T,H)^T$ at $z_0$,
where $T,H\in\R$, and set $s=z_0-z_{-1}$ and $y=g-z_{-1}$.
Then $g^Ts=0$, $s^Ty=L^2$, and
\begin{equation}\label{eq:M}
\U(D;s,y)=\mathcal M_{L,a,b}(T,H)
:=\diag(0,a,b)+\frac1{L^2}(-L,T,H)^T(-L,T,H)\succ0.
\end{equation}
If $R+T>0$, then $(z_0,g,\mathcal M_{L,a,b}(T,H))\in\mathcal A$.
Also, $\mathcal M_{L,a,b}(0,0)=\diag(1,a,b)$.
\end{lemma}

\begin{proof}
The converse in Lemma~\ref{lem:canonical} gives $D\succ0$,
$\lambda_Q(z_{-1},D)=\theta$ and $Z_Q(z_{-1},D)=z_0$.
The update calculation in Appendix~\ref{app:incoming} gives \eqref{eq:M} and
$z_0^T\mathcal M_{L,a,b}(T,H)^{-1}g=R(R+T)/a$,
proving the remaining assertions.
\end{proof}

\subsection{Constructing one cycle}\label{sec:cycle}\label{sec:local}

Each cycle produces a large quotient and prepares the next cycle.
A $D$-invariant plane is a two-dimensional linear subspace
$P\subset\R^3$ satisfying $Dv\in P$ for every $v\in P$. Set

\begin{equation}\label{eq:rho}
\rho(b):=\frac{b-1}{\sqrt{2(b^2+1)}}\qquad(b>1).
\end{equation}

\begin{definition}[Starting data for cycle $m$]\label{def:cycle-input}
Let $m\ge1$ be an integer. A point $z_{-1}$ with prescribed gradient
$z_{-1}$ and matrix $D$ is admissible for cycle $m$ if:
\begin{enumerate}
\item $z_{-1}\ne0$ and $D=D^T\succ0$.
\item There is a two-dimensional $D$-invariant subspace
$P\subset\R^3$ containing $z_{-1}$. Since $D$ is symmetric,
$P^\perp$ is a one-dimensional invariant subspace; denote the
eigenvalue of $D$ on it by $b$.
\item The point $z_0:=Z_Q(z_{-1},D)$ and the eigenvalue $b$ satisfy
\begin{equation}\label{eq:cycle-state}
1<b<11/10,\qquad
0<\norm{z_0}<\norm{z_{-1}}/8,\qquad
\rho(b)/\norm{z_0}\ge4m.
\end{equation}
\end{enumerate}
\end{definition}

The one-step properties and form \eqref{eq:canonical} follow from
Lemmas~\ref{lem:transition} and~\ref{lem:canonical}.
Table~\ref{tab:cycle} summarizes the cycle: the gradients are
changed first at $z_0$, then at $z_3$, and the output
$(z_4,z_4,D_{\rm new})$ starts the next cycle.

\begin{table}[!ht]
\caption{Point, prescribed gradient and matrix data in one cycle.
Only the gradients at $z_0$ and $z_3$ are changed from those of $Q$.}
\label{tab:cycle}\label{fig:proof-architecture}
\centering
\small
\setlength{\tabcolsep}{4pt}
\renewcommand{\arraystretch}{1.12}
\begin{tabularx}{\textwidth}{@{}>{\raggedright\arraybackslash}p{.17\textwidth}
>{\raggedright\arraybackslash}p{.38\textwidth}
>{\raggedright\arraybackslash}X@{}}
\toprule
Step & Input $\longrightarrow$ output data & Role\\
\midrule
$z_{-1}\to z_0$
& $(z_{-1},z_{-1},D)\to(z_0,g_0,M)$
& Prescribe $g_0$ and compute $M$ by the incoming BFGS update (see Lemma~\ref{lem:bridge}).\\
\addlinespace
$z_0\to z_1$
& $(z_0,g_0,M)\to(z_1,z_1,M_1)$
& Begin the three-step quotient construction.\\
\addlinespace
$z_1\to z_2$
& $(z_1,z_1,M_1)\to(z_2,z_2,M_2)$
& Use the quadratic gradient at $z_1$.\\
\addlinespace
$z_2\to z_3$
& $(z_2,z_2,M_2)\to(z_3,g_3,B_{\rm reg})$
& Fix the quotient; choose gradient $g_3\ne z_3$ to avoid a zero next point and recompute the matrix.\\
\addlinespace
$z_3\to z_4$
& $(z_3,g_3,B_{\rm reg})\to(z_4,z_4,D_{\rm new})$
& Obtain starting data for the next cycle.\\
\bottomrule
\end{tabularx}
\end{table}

\FloatBarrier
\subsubsection{The three-step quotient}\label{sec:bad}

With $g_0=z_0$, the next point would be $0$. We perturb this
gradient in two directions: a change of order $\tau$ in the
invariant plane and a change of order $\tau^3$ normal to it.

\begin{lemma}[Limit of the three-step quotient]\label{lem:bad}

Fix $L,R,a>0$ and $b>1$. For each $\tau>0$, set
\[
\kappa:=-\frac{b^2R^3}{a(b-1)L^2},\qquad
T:=\tau R,\qquad H:=\kappa\tau^3.
\]
Let $z_0:=(0,R,0)^T$, $g_0:=(0,R+T,H)^T$, and
$M:=\mathcal M_{L,a,b}(T,H)$, with $\mathcal M_{L,a,b}$ defined
in \eqref{eq:M}. Using the map $\Sstep$ in
\eqref{eq:transition}, define the three-step sequence by
\[
(z_1,M_1):=\Sstep(z_0,g_0,M),\qquad
(z_{i+1},M_{i+1}):=\Sstep(z_i,z_i,M_i),\quad i=1,2.
\]
There exists $\tau_0=\tau_0(L,R,a,b)>0$ such that all three steps
are well defined and $z_1,z_2,z_3\ne0$ for $0<\tau<\tau_0$.
Moreover, with $\rho$ defined in \eqref{eq:rho},
\begin{equation}\label{eq:quotient}
\frac{\norm{z_3}}{\norm{z_2}\norm{z_0}}
\longrightarrow\frac{\rho(b)}R
\qquad\text{as }\tau\downarrow0.
\end{equation}
\end{lemma}

The proof is given in Appendix~\ref{app:normal}. For fixed starting
data, reducing $\tau$ approaches the limit $\rho(b)/R$; to make
this limit larger in later cycles, we must also improve the
starting data.

\Needspace{16\baselineskip}
\subsubsection{Preparing the next cycle}\label{sec:regeneration}\label{sec:interface}

The preceding lemma temporarily uses gradient $z_3$ at the last
point. Setting $g_3=z_3$ would make the next point $0$. We instead
choose $g_3$ perpendicular to the fixed displacement $z_3-z_2$
and recompute the matrix at $z_3$. This leaves the current
quotient fixed and allows one further step to prepare the next cycle.

\begin{lemma}[Starting the next cycle]\label{lem:bridge}\label{lem:regeneration}
Fix an integer $m\ge1$ and $L,R,a>0$, $b>1$.
For every sufficiently small $\tau>0$, take the reference data
$z_1,z_2,z_3,M_2$ from Lemma~\ref{lem:bad}.
With these data fixed, for every $\eta_*>0$ one can choose a gradient
$g_3$ at $z_3$ with $\norm{g_3-z_3}\le\eta_*$, preserving
$g_3^T(z_3-z_2)=0$, such that the actual incoming update and an
admissible outgoing step
\[
B_{\rm reg}:=\U(M_2;z_3-z_2,g_3-z_2),\qquad
(z_4,D_{\rm new}):=\Sstep(z_3,g_3,B_{\rm reg}),
\]
produce starting data $(z_4,z_4,D_{\rm new})$ for the next cycle,
satisfying Definition~\ref{def:cycle-input} for cycle $m+1$
and $0<\norm{z_4}<\norm{z_3}/8$.
\end{lemma}

The key is to choose the direction of the gradient change so that
the actual updated matrix has an invariant plane containing
the new point. This gives a scalar equation with a simple
limiting root, which persists by the implicit function theorem.
Decreasing the perturbation then makes the gradient error and
the new point arbitrarily small while increasing the next
cycle's limiting quotient. The proof is in
Appendix~\ref{app:variation}. In Section~\ref{sec:objective}, we will combine these two lemmas to prove
Lemma~\ref{lem:trajectory}.

\section{An infinite trajectory and a smooth objective}\label{sec:objective}

We first join the cycles from Section~\ref{sec:cycle} into an
infinite sequence to prove Lemma~\ref{lem:trajectory}. We then
construct a smooth objective whose BFGS iterates follow this
sequence to prove Lemma~\ref{lem:realization}.

\subsection{Proof of Lemma~\ref{lem:trajectory}}\label{sec:trajectory}

The local results provide a repeatable cycle for each fixed input.
We impose the global
requirements from Section~\ref{sec:realization}: the indexed
quotient bounds, separation, and errors small enough for smooth
realization. The tolerance constant $c>0$ is fixed before the
recursion; each subsequent choice concerns only the current
finite data.

\begin{proof}
Fix $c>0$ throughout the construction. Take \eqref{eq:canonical} with
\[
L:=\frac{1}{4},\qquad R:=\frac{1}{10000},\qquad
\theta:=a:=1,\qquad b:=\frac{101}{100}.
\]
Lemma~\ref{lem:canonical} gives $D\succ0$ and
$z_0=Z_Q(z_{-1},D)=(0,R,0)^T$, with invariant plane $\Span(e_1,e_2)$.
Since $\norm{z_{-1}}<1/2$ and
\[
\frac{\norm{z_0}}{\norm{z_{-1}}}\le\frac RL=\frac1{2500}<\frac18,
\qquad
\frac{\rho(b)}{\norm{z_0}}=\frac{10000}{\sqrt{40402}}>4,
\]
the triple $(x_1,g_1,B_1):=(z_{-1},z_{-1},D)$ satisfies
Definition~\ref{def:cycle-input} for $m=1$. Set $O_1:=I$.

Suppose the triples through index $5m-4$ are fixed, with
$g_{5m-4}=x_{5m-4}$, and the last triple satisfies
Definition~\ref{def:cycle-input} for cycle $m$.
In an adapted frame $O_m$ from Lemma~\ref{lem:canonical}, write
\[
z_{-1}:=O_m^Tx_{5m-4},\qquad D:=O_m^TB_{5m-4}O_m.
\]
These data have form \eqref{eq:canonical}, with fixed
$L,R,\theta,a>0$, $1<b<11/10$, and $z_0=(0,R,0)^T$.
Local data below use this frame; the indexed triples remain in
fixed ambient coordinates.

Choose $\tau=\tau_m>0$ in the construction of Lemma~\ref{lem:bad}
so that both Lemmas~\ref{lem:bad} and \ref{lem:bridge} apply and
\begin{equation}\label{eq:tau-choice}
\begin{gathered}
\norm{g_0-z_0}\le c\,2^{-(5m-3)}R^{5m},\qquad
\frac{\norm{z_3}}{\norm{z_2}R}
 \ge\frac{\rho(b)}{2R}\ge2m,\\
\norm{z_1}<\frac{R}{8},\qquad
\norm{z_2}<\frac{1}{8}\norm{z_1},\qquad
\norm{z_3}<\frac{1}{8}\norm{z_2}.
\end{gathered}
\end{equation}
For this fixed predecessor, $\norm{g_0-z_0}\to0$ and
\eqref{eq:local-ratios} gives the first two norm ratios tending
to $0$ and the last to $\rho(b)<(b-1)/2<1/20$.
The quotient limit \eqref{eq:quotient} is $\rho(b)/R\ge4m$.
Thus these bounds and the two lemmas' smallness requirements
hold simultaneously on a sufficiently small positive $\tau$-interval.

Fixing $\tau_m$ determines $g_0,M,z_1,z_2,z_3,M_1,M_2$;
the gradient and incoming matrix at $z_3$ remain unassigned.
Apply Lemma~\ref{lem:bridge} with the now fixed tolerance
$c\,2^{-5m}\norm{z_3}^{5m+3}>0$, choosing its parameter
$\epsilon_m$ after $\tau_m$. It gives $g_3,B_{\rm reg},z_4,D_{\rm new}$, with
$(z_4,z_4,D_{\rm new})$ satisfying Definition~\ref{def:cycle-input}
for cycle $m+1$, and
\begin{equation}\label{eq:epsilon-choice}
\norm{g_3-z_3}\le c\,2^{-5m}\norm{z_3}^{5m+3},\qquad
0<\norm{z_4}<\norm{z_3}/8.
\end{equation}
It preserves the preceding data and quotient, and verifies
\begin{equation}\label{eq:changed-endpoint}
\begin{gathered}
g_3^T(z_3-z_2)=0,\qquad
(z_3-z_2)^T(g_3-z_2)=\norm{z_3-z_2}^2>0,\\
B_{\rm reg}=\U(M_2;z_3-z_2,g_3-z_2).
\end{gathered}
\end{equation}
Assign the five new triples in ambient coordinates by
\begin{equation}\label{eq:indices}
\begin{aligned}
(x_{5m-3},g_{5m-3},B_{5m-3})
 &:= (O_mz_0,O_mg_0,O_mMO_m^T),\\
(x_{5m-2},g_{5m-2},B_{5m-2})
 &:= (O_mz_1,O_mz_1,O_mM_1O_m^T),\\
(x_{5m-1},g_{5m-1},B_{5m-1})
 &:= (O_mz_2,O_mz_2,O_mM_2O_m^T),\\
(x_{5m},g_{5m},B_{5m})
 &:= (O_mz_3,O_mg_3,O_mB_{\rm reg}O_m^T),\\
(x_{5m+1},g_{5m+1},B_{5m+1})
 &:= (O_mz_4,O_mz_4,O_mD_{\rm new}O_m^T).
\end{aligned}
\end{equation}
In particular, the matrix at $z_3$ is the actual $B_{\rm reg}$,
not the auxiliary $M_3$.

The incoming step $z_{-1}\to z_0$ produces $M$ by
Lemma~\ref{lem:endpoint-update}; the first two calls in
Lemma~\ref{lem:bad} give $z_0\to z_1\to z_2$; and
Lemma~\ref{lem:bridge} gives the last two actual steps, with
\eqref{eq:changed-endpoint} preserving the displacement
$z_2\to z_3$. Their positive line parameters are
\[
\beta_j:=\frac{x_j^TB_j^{-1}g_j}{\norm{B_j^{-1}g_j}^2}
\qquad(5m-4\le j\le5m),
\]
with $\beta_{5m-4}=\theta$. By Lemma~\ref{lem:transition},
orthogonal transport in \eqref{eq:indices} therefore gives
\eqref{eq:algebraic-run}, including positive curvature and
positive definite matrices. All points are nonzero by the local
lemmas. All gradients are nonzero as well: $g_0$ has second
component $R(1+\tau_m)>0$, $g_3$ belongs to an admissible outgoing
call, and the others equal their points.

Only indices $5m-3$ and $5m$ can have gradient errors.
Orthogonality and \eqref{eq:tau-choice}--\eqref{eq:epsilon-choice}
give the error bound in \eqref{eq:flat-data}, since the local
radii are $R$ and $\norm{z_3}$ and the exponents $5m$ and $5m+3$
are exactly $j+3$ at these indices.
The incoming bound \eqref{eq:cycle-state}, the three block bounds
\eqref{eq:tau-choice}, and \eqref{eq:epsilon-choice} cover all five
contractions in \eqref{eq:flat-data}. The quotient bound in
\eqref{eq:tau-choice} becomes \eqref{eq:bad-subsequence} under
\eqref{eq:indices}.

Since $5(m+1)-4=5m+1$, the last triple is the next cycle's
starting data. Set $O_{m+1}:=O_mO_{\rm new}$ using
\eqref{eq:bridge-frame}--\eqref{eq:bridge-canonical}.
Its next endpoint $x_{5m+2}=Z_Q(x_{5m+1},B_{5m+1})$ is fixed,
but its gradient and incoming matrix are chosen only with
$\tau_{m+1}$; the actual incoming update forms the direct input
of Lemma~\ref{lem:bad} as verified in \eqref{eq:bridge-input}.
Thus each stage extends the fixed prefix by five triples and
five positive step parameters without changing earlier assignments.
The choices occur in the order $\tau_m,\epsilon_m,\tau_{m+1}$,
with $c$ fixed from the start. All local thresholds may depend
on already fixed data, so induction defines one infinite sequence
without a uniform threshold over cycles. Finally,
\eqref{eq:flat-data} implies
$0<\norm{x_j}<\tfrac12\,8^{-(j-1)}$, proving distinctness and
convergence to $0$.
\end{proof}

\subsection{Proof of Lemma~\ref{lem:realization}}\label{sec:smooth-objective}\label{sec:identification}

We first construct the function with the prescribed gradients,
then show that its BFGS iteration follows the prescribed sequence.

\begin{proof}
Fix a smooth cutoff $\chi:\R^3\to[0,1]$ that equals one on
$\{\xi:\norm{\xi}\le1/2\}$ and vanishes on
$\{\xi:\norm{\xi}\ge1\}$; an explicit choice is given in
Appendix~\ref{app:interpolation}. For each unit vector $u$, put
$f_u(\xi):=(u^T\xi)\chi(\xi)$ and define
\[
K_2:=\sup_{\substack{u,\xi\in\R^3\\\norm{u}=1}}
       \norm{\nabla^2f_u(\xi)}<\infty,
\qquad c_0:=\frac1{8(1+K_2)}.
\]
The constant $c_0$ depends only on the fixed cutoff. Take any data
as in Lemma~\ref{lem:trajectory} with $c=c_0$. Write
\[
r_j:=\norm{x_j},\qquad h_j:=\frac{r_j}{4},\qquad \delta_j:=g_j-x_j.
\]
Define the local gradient corrections and the objective by
\begin{equation}\label{eq:bumps}
\begin{aligned}
\psi_j(x)&:=\delta_j^T(x-x_j)\chi\bigl((x-x_j)/h_j\bigr),\\
F(x)&:=Q(x)+\sum_{j=1}^\infty\psi_j(x).
\end{aligned}
\end{equation}
Each correction is supported in the closed ball of radius $h_j$
about $x_j$. Appendix~\ref{app:interpolation} proves that these
balls are pairwise disjoint, accumulate only at $0$, and have
closure inside the ball of radius $5r_1/4$. It also establishes
\[
F(x_j)=Q(x_j),\qquad \nabla F(x_j)=g_j,
\]
and, using \eqref{eq:bump-bound} with $\ell=2$, we have 
$\norm{\nabla^2(F-Q)}_\infty\le4K_2c_0<1/2$ away from $0$.
The only remaining analytic issue is smoothness at the accumulation
point. Appendix~\ref{app:interpolation} proves that every derivative
of these corrections tends to $0$ faster than any fixed power of
$\norm{x}$, using the increasing exponent in \eqref{eq:flat-data}.
Consequently, $F-Q$ is smooth, compactly supported and flat at $0$,
and the Hessian bounds hold everywhere.

Flatness gives $F(0)=0$, $\nabla F(0)=0$ and $\nabla^2F(0)=I$.
For every $x\in\R^3$,
\[
\begin{aligned}
x^T\nabla F(x)
&=\int_0^1x^T\nabla^2F(tx)x\,dt
\ge\tfrac12\norm{x}^2,\\
F(x)-F(0)
&=\int_0^1(1-t)x^T\nabla^2F(tx)x\,dt
\ge\tfrac14\norm{x}^2.
\end{aligned}
\]
Hence $0$ is the unique minimizer and no other gradient vanishes.

It remains to identify the BFGS iterates. The interpolation in
\eqref{eq:bumps} supplies the prescribed gradient at every point
$x_j$, and strong convexity turns line-search stationarity into
unique minimization. We proceed by induction. At a prescribed point $x_j$,
$d_j=-B_j^{-1}g_j\ne0$, and \eqref{eq:algebraic-run} gives
\[
x_{j+1}=x_j+\beta_jd_j,\qquad
\nabla F(x_{j+1})^Td_j
=\beta_j^{-1}g_{j+1}^Ts_j=0.
\]
For $\varphi_j(\alpha)=F(x_j+\alpha d_j)$, the Hessian bound gives
\[
\varphi_j'(\beta_j)=0,\qquad
\varphi_j''(\alpha)\ge\tfrac12\norm{d_j}^2>0.
\]
Integrating on either side of $\beta_j$ yields
\[
\varphi_j(\alpha)\ge\varphi_j(\beta_j)
+\tfrac14\norm{d_j}^2(\alpha-\beta_j)^2
\qquad(\alpha\in\R).
\]
Thus $\beta_j>0$ is the unique minimizing line-search parameter.
The new point is $x_{j+1}$, its gradient is $g_{j+1}$, and the
gradient difference is $y_j$,
so the next actual matrix is the prescribed $B_{j+1}$ by
\eqref{eq:algebraic-run}. Starting with the common $(x_1,B_1)$ proves
the identification for every $j$. The quotient bound
\eqref{eq:bad-subsequence} is therefore preserved.
The positive curvature products and positive definite matrices
are therefore those of the actual iteration. Also, directly from
the Hessian bound,
\[
s_j^Ty_j
=\int_0^1s_j^T\nabla^2F(x_j+ts_j)s_j\,dt
\ge\tfrac12\norm{s_j}^2>0.
\]
The first bound in \eqref{eq:flat-data} gives nonzero, distinct
iterates converging to $0$.
\end{proof}

\section{Conclusion}\label{sec:conclusion}

We have disproved the three-step product estimate \eqref{eq:product}
for a globally strongly convex $C^\infty$ objective on $\R^3$ that
agrees with $Q$ to every order at its minimizer. The five-step
construction and its smooth realization yield a single convergent
BFGS trajectory with exact line searches whose product quotients
have infinite limsup.

\section*{Acknowledgment}

Generative artificial intelligence tools were used
to support the development, checking, and
exposition of mathematical proofs. The authors critically
assessed and revised the arguments and retain full responsibility
for the mathematical content and the final manuscript.

\bibliographystyle{siamplain}
\bibliography{final_draft}

\newpage

\appendix
\section{Proof of Lemma~\ref{lem:endpoint-update}}\label{app:incoming}

We first justify the coordinate form \eqref{eq:canonical} for
the invariant plane used in the construction, then give the
matrix calculation for Lemma~\ref{lem:endpoint-update}.

\begin{lemma}[Coordinates for the incoming step]\label{lem:canonical}
Let $D=D^T\succ0$, let $z_{-1}\ne0$ belong to a two-dimensional $D$-invariant
linear subspace $P\subset\R^3$, and suppose $z_{-1}$ is not an eigenvector of $D$.
Set
\[
z_0:=Z_Q(z_{-1},D),\qquad L:=\norm{z_{-1}-z_0},\qquad R:=\norm{z_0},
\qquad \theta:=\lambda_Q(z_{-1},D).
\]
The step length $L$, point norm $R$ and line parameter $\theta$
are positive.
In the orthonormal frame with
first two axes $(z_{-1}-z_0)/L$ and $z_0/R$, and third axis a unit normal to
$P$, let $b>0$ denote the normal eigenvalue of $D$.
For some $a>0$, the data have the form \eqref{eq:canonical}.
Conversely, any $L,R,\theta,a,b>0$ in \eqref{eq:canonical} give
$D\succ0$, $\lambda_Q(z_{-1},D)=\theta$ and $Z_Q(z_{-1},D)=z_0$.

\end{lemma}

\begin{proof}[Proof of Lemma~\ref{lem:canonical}]
The non-eigenvector hypothesis gives $z_0\ne0$ by
\eqref{eq:quadratic-endpoint}. Also $z_{-1}-z_0=\theta D^{-1}z_{-1}\ne0$ and
$\theta>0$. Since the restriction of $D$ to $P$ is invertible,
$P$ is invariant under $D^{-1}$, so both $z_{-1}-z_0$ and $z_0$ belong to
$P$. The identity $z_0^T(z_0-z_{-1})=0$ makes their unit directions orthogonal.
Symmetry makes $P^\perp$ invariant as well, with a positive
eigenvalue $b$. In this frame,
\[
De_1=(\theta/L)z_{-1}=\theta e_1+(\theta R/L)e_2.
\]
Symmetry determines the first row. The remaining coefficient is
\[
a:=e_2^TDe_2-\theta R^2/L^2
=\bigl(e_2-(R/L)e_1\bigr)^TD\bigl(e_2-(R/L)e_1\bigr)>0.
\]
This proves \eqref{eq:canonical}. Conversely, that expression gives
\[
v^TDv=\theta\bigl(v_1+(R/L)v_2\bigr)^2+av_2^2+bv_3^2>0
\quad(v\ne0),
\qquad D^{-1}z_{-1}=(L/\theta)e_1.
\]
Substitution in \eqref{eq:quadratic-endpoint} gives the stated line
parameter and next point.

\end{proof}

\begin{proof}[Matrix calculation for Lemma~\ref{lem:endpoint-update}]
In these coordinates
\[
s=(-L,0,0)^T,\quad y=(-L,T,H)^T,\quad
Ds=-\theta z_{-1},\quad s^TDs=\theta L^2,\quad s^Ty=L^2.
\]
Thus the first rank-one term in the BFGS update cancels exactly
$\theta z_{-1}z_{-1}^T/L^2$ in $D$. This proves \eqref{eq:M} and explains
why the outgoing matrix does not depend on $\theta$.
Its positive definiteness follows from Lemma~\ref{lem:transition}.
For $w:=\mathcal M_{L,a,b}(T,H)^{-1}g$, the component formulas are
\begin{equation}\label{eq:inverse}
w_2=\frac{R+T}{a},\qquad w_3=\frac H b,\qquad
w_1=\frac{Tw_2+Hw_3}{L}.
\end{equation}
The first row of $\mathcal M w=g$ gives
$w_1=(Tw_2+Hw_3)/L$; the other two then give
$aw_2=R+T$ and $bw_3=H$. This proves \eqref{eq:inverse} and the
remaining assertions, since $z_0^Tw=R(R+T)/a$.
\end{proof}

\section{Proof of Lemma~\ref{lem:bad}}\label{app:normal}

We first establish the normal perturbation estimate used in the
three-step construction. The second line-search denominator
vanishes at the limiting data, so an explicit bound is needed.

\begin{lemma}\label{lem:normal-response}
Fix $L,R,a>0$, $b>1$, and $T=\tau R$ as in Lemma~\ref{lem:bad},
but allow $H$ to vary independently. For some $\tau_0>0$, the first
two calls to $\Sstep$ are well defined for
$0<\tau<\tau_0$, $|H|\le\tau$.
The first output and matrix are analytic near $(\tau,H)=(0,0)$,
with $M_1(0,0)=\diag(1,1,b)$ and
\[
z_{1,3}=-aH/b+O(|\tau H|+|H|^3).
\]
Writing $z_2^0=z_2(\tau,0)$, the second output satisfies
\eqref{eq:z2-perturbation}, with constants depending only on the
fixed parameters. In particular, its planar error is $O(H^2/\tau)$
and its leading normal response is $-a(b-1)H/b^2$.
\end{lemma}

\begin{proof}
Regard $\tau$ and $H$ as independent variables, retain
$T=\tau R$, and denote the first transition's outputs by
$z_1(\tau,H)$ and $M_1(\tau,H)$.

$w_0$ and $\lambda_0$ denote the direction vector and line parameter
of the first transition, now depending on $(\tau,H)$. Their planar
values are given by \eqref{eq:planar}.

 At $(\tau,H)=(0,0)$, one has
\[
w_0=(R/a)e_2,\qquad \lambda_0=a,
\qquad s_0:=z_1-z_0=-Re_2,
\qquad y_0:=z_1-g_0=-Re_2.
\]
The three denominators $w_0^Tw_0$, $s_0^T\mathcal M_{L,a,b}(0,0)s_0$, and $s_0^Ty_0$
are respectively $R^2/a^2$, $aR^2$, and $R^2$. Hence the rational
formulas defining $z_1$ and $M_1$ are analytic in a neighborhood of
$(0,0)$, where
\[
z_1(0,0)=0,
\qquad M_1(0,0)=\diag(1,1,b)\succ0.
\]
In particular, after shrinking the neighborhood, $M_1$ and $M_1^{-1}$
are uniformly bounded and $M_1$ is uniformly positive definite. The value $z_1(0,0)=0$ is used only for these expansions.

Reflection by $J:=\diag(1,1,-1)$ sends $H$ to $-H$, so
$z_1(\tau,-H)=Jz_1(\tau,H)$ and
$M_1(\tau,-H)=JM_1(\tau,H)J$. Write a superscript $0$ for evaluation
at $H=0$. Before using block notation, for a vector $r$ and a matrix
$G$ define
\[
r_{1:2}:=(r_1,r_2)^T,
\quad G_{1:2,1:2}:=(G_{ij})_{i,j=1}^2,
\quad G_{1:2,3}:=(G_{13},G_{23})^T.
\]
Analyticity and parity give, uniformly near $(0,0)$,
\begin{equation}\label{eq:parity}
\begin{aligned}
(z_1)_{1:2}-(z_1^0)_{1:2}&=O(H^2),&
z_{1,3}&=-aH/b+O(|\tau H|+|H|^3),\\
(M_1)_{1:2,1:2}-(M_1^0)_{1:2,1:2}&=O(H^2),&
(M_1)_{1:2,3}&=O(|H|),\\
(M_1)_{33}&=b+O(H^2).
\end{aligned}
\end{equation}
Here $z_{1,3}$ denotes the third component. Its formula uses $z_{1,3}=-\lambda_0H/b$ and
$\lambda_0=a+O(|\tau|+H^2)$. The inverse $M_1^{-1}$ has the same
parity and block estimates.

For the inverse, ``the same'' refers to parity and orders of error;
its leading normal entry is $b^{-1}$, rather than $b$.

At the limiting parameter value, $z_1=0$ makes the second transition degenerate, so the analyticity
of the first transition cannot simply be carried over. We retain explicitly the loss of scale caused
by the vanishing denominator.

Restrict to $0<\tau<\tau_0$ and $|H|\le\tau$, with $\tau_0>0$
sufficiently small. In what follows, $f\asymp\tau$ means that
$c\tau\le f\le C\tau$ for constants $c,C>0$ independent of
$(\tau,H)$ in this region. Put $w_1:=M_1^{-1}z_1$. The first component
in \eqref{eq:planar}, its $O(H^2)$ perturbation, and the uniform bounds
on $M_1$ and $M_1^{-1}$ imply
\[
\norm{z_1}\asymp\tau,
\qquad \norm{w_1}\asymp\tau.
\]
Moreover,
\[
(w_1)_{1:2}-(w_1^0)_{1:2}=O(H^2),\qquad
w_{1,3}=z_{1,3}/b+O(|\tau H|+|H|^3).
\]
Define
\[
n_1:=z_1^Tw_1,
\qquad d_1:=w_1^Tw_1,
\qquad \lambda_1:=n_1/d_1,
\]
and use a superscript $0$ for their values at $H=0$. The preceding
estimates give
\[
n_1-n_1^0=O(H^2),\quad d_1-d_1^0=O(H^2),\quad
n_1^0=O(\tau^2),\quad d_1,d_1^0\asymp\tau^2.
\]
Consequently,
\begin{equation}\label{eq:lambda-control}
\lambda_1-\lambda_1^0
=\frac{(n_1-n_1^0)d_1^0-n_1^0(d_1-d_1^0)}{d_1d_1^0}
=O(H^2/\tau^2).
\end{equation}
At $H=0$, the analytic vectors $z_1^0/\tau$ and $w_1^0/\tau$ both
tend to $-(R^2/L)e_1$, so $\lambda_1^0=1+O(\tau)$. Hence
\[
\lambda_1=1+O\!\left(\tau+\frac{H^2}{\tau^2}\right).
\]
Using $z_2=z_1-\lambda_1w_1$ gives
\begin{equation}\label{eq:z2-perturbation}
\begin{aligned}
(z_2)_{1:2}-(z_2^0)_{1:2}&=O(H^2/\tau),\\
z_{2,3}
&=\left(1-\frac{\lambda_1}{b}\right)z_{1,3}
  +O(|\tau H|+|H|^3)\\
&=-\frac{a(b-1)}{b^2}H
  +O\!\left(\tau|H|+\frac{|H|^3}{\tau^2}\right).
\end{aligned}
\end{equation}

The bounds $\norm{z_1}\asymp\tau$ and $M_1\succ0$ ensure that the
second call is admissible and its projection denominator is
positive. The first is admissible because $R+\tau R>0$.
Lemma~\ref{lem:transition} supplies their update denominators and
positive curvature. This proves the stated domain as well as the
estimates.
\end{proof}

\begin{proof}[Proof of Lemma~\ref{lem:bad}]

All estimates below keep $L,R,a,b$ fixed, and all limits are taken as
$\tau\downarrow0$. Put $A:=R^3/L^2$. We first compute the leading
vectors and the second updated matrix, then take norms to obtain
the product quotient.

First consider the auxiliary planar family obtained
by setting $H=0$ while retaining $T=\tau R$. Put
\[
\eta:=R/L,\qquad c_\tau:=1+\tau,
\qquad d:=1+\tau^2\eta^2.
\]
With $g_0=c_\tau z_0$, define
$w_0:=\mathcal M_{L,a,b}(T,0)^{-1}g_0$ and
$\lambda_0:=z_0^Tw_0/\norm{w_0}^2$. Formula \eqref{eq:inverse} and
\eqref{eq:transition} give
\begin{equation}\label{eq:planar}
w_0=\frac{c_\tau R}{a}(\tau\eta,1,0)^T,\qquad
\lambda_0=\frac{a}{c_\tau d},\qquad
z_1=\frac R d(-\tau\eta,\tau^2\eta^2,0)^T.
\end{equation}
Let $y_0:=z_1-g_0=z_1-c_\tau z_0$ and
$s_1:=z_2-z_1$. The first transition has positive curvature because
$R+T>0$. Formula \eqref{eq:planar} gives $z_1\ne0$, so the second
transition is well defined. All relevant vectors lie in $\Span(e_1,e_2)$.
To determine $z_2$ without expanding a second inverse matrix, we
use the following local consequence of the BFGS secant equation.
For any adjacent steps, suppose $B_j=B_j^T\succ0$,
$s_j^Ty_j>0$, and $g_{j+1}^Ts_j=0$. Set
$B_{j+1}=\U(B_j;s_j,y_j)$ and
$s_{j+1}=-\alpha_{j+1}B_{j+1}^{-1}g_{j+1}$, where
$\alpha_{j+1}\in\R$. The standard BFGS properties give
$B_{j+1}\succ0$ and $B_{j+1}s_j=y_j$, hence
\begin{equation}\label{eq:conjugacy}
s_{j+1}^Ty_j
=-\alpha_{j+1}g_{j+1}^TB_{j+1}^{-1}y_j
=-\alpha_{j+1}g_{j+1}^Ts_j=0.
\end{equation}
This assertion is independent of the indexing origin; applying it
to the step from $z_0$ to $z_1$, with $B_0=M$, $B_1=M_1$,
$s_0=z_1-z_0$, and $g_1=z_1$, gives $s_1\perp y_0$.
Also $z_2\perp s_1$ by \eqref{eq:compatibility}. Since $s_1\ne0$ and $y_0\ne0$, it follows
that $z_2$ is the orthogonal projection of $z_1$ onto
$\Span(y_0)$. The identity $z_1^T(z_1-z_0)=0$ gives
$z_1^Tz_0=\norm{z_1}^2$, and hence
\begin{equation}\label{eq:planar-z2}
z_2=-\tau\frac{\norm{z_1}^2}{\norm{y_0}^2}y_0,\qquad
\norm{z_1}^2=\frac{R^2\tau^2\eta^2}{d},\qquad
y_0\to-Re_2.
\end{equation}
Thus $\tau^{-3}z_2\to Ae_2$ in the planar family.

Lemma~\ref{lem:normal-response} controls the two-parameter family while retaining the vanishing
denominator of the second transition. For a vector $r$, write
$r_{1:2}=(r_1,r_2)^T$ for its first two components and $r_3$ for
its third component. For $0<\tau<\tau_0$ and $|H|\le\tau$, the lemma gives the
estimates \eqref{eq:z2-perturbation}.
The constants depend on fixed $L,R,a,b$. The planar error is
smaller than order $\tau^3$ when $H=\kappa\tau^3$.
The choice of $\kappa$ makes the normal linear response equal to
the planar leading term $A\tau^3$.

For the prescribed curve $H=\kappa\tau^3$, the restriction
$|H|\le\tau$ holds for all sufficiently small $\tau>0$. Along this
curve, the first error in \eqref{eq:z2-perturbation} is $O(\tau^5)$,
and the definition of $\kappa$ gives
$z_{2,3}=A\tau^3+O(\tau^4)$.
The first-output bounds in \eqref{eq:parity} also give
$z_1(\tau,\kappa\tau^3)=z_1(\tau,0)+O(\tau^3)$.
Analyticity of the first updated matrix gives
$M_1\to\diag(1,1,b)$ along this curve.
Together with
\eqref{eq:planar}--\eqref{eq:planar-z2}, these estimates give
\begin{equation}\label{eq:z1}
\tau^{-1}z_1\longrightarrow-\frac{R^2}{L}e_1
\end{equation}
and
\begin{equation}\label{eq:z2}
\tau^{-3}z_2\longrightarrow A(e_2+e_3).
\end{equation}

For the second transition, the gradient difference is
$y_1:=z_2-z_1=s_1$. Define $h_1:=s_1/\norm{s_1}$. The established
limits give $h_1\to e_1$, and therefore
\begin{equation}\label{eq:M2-limit}
M_2=M_1-\frac{M_1h_1h_1^TM_1}{h_1^TM_1h_1}+h_1h_1^T
\longrightarrow\diag(1,1,b).
\end{equation}
Indeed, $h_1^TM_1h_1\to1$ and both rank-one terms tend to
$e_1e_1^T$. The nonzero limits in \eqref{eq:z1} and \eqref{eq:z2}
ensure that $z_1,z_2\ne0$ for all sufficiently small $\tau>0$.

For the third transition, define
\[
w_2:=M_2^{-1}z_2,
\qquad \lambda_2:=\frac{z_2^Tw_2}{\norm{w_2}^2}.
\]
Then
\[
\tau^{-3}w_2\to A(e_2+b^{-1}e_3),
\qquad
\tau^{-6}\norm{w_2}^2\to A^2(1+b^{-2})>0,
\qquad
\lambda_2\to\frac{b(b+1)}{b^2+1}.
\]
Thus $z_3=z_2-\lambda_2w_2$ gives
\begin{equation}\label{eq:z3}
\tau^{-3}z_3\longrightarrow
A\frac{b-1}{b^2+1}(-e_2+be_3).
\end{equation}
Its nonzero limit gives $z_3\ne0$ for all sufficiently small $\tau>0$.
Lemma~\ref{lem:transition} now gives positive curvature for each of the
three transitions, positive line parameters and positive definite
updated matrices.
Taking norms in \eqref{eq:z2} and \eqref{eq:z3}, and using
$\norm{z_0}=R$, proves \eqref{eq:quotient}.
The same limits, together with \eqref{eq:z1}, give
\begin{equation}\label{eq:local-ratios}
\frac{\norm{z_1}}{\norm{z_0}}\longrightarrow0,\qquad
\frac{\norm{z_2}}{\norm{z_1}}\longrightarrow0,\qquad
\frac{\norm{z_3}}{\norm{z_2}}\longrightarrow\rho(b).
\end{equation}
\end{proof}

\Needspace{6\baselineskip}
\section{Proof of Lemma~\ref{lem:bridge}}\label{app:variation}

\begin{proof}[Proof of Lemma~\ref{lem:bridge}]
Keep the original $L,R,a>0$ and $b>1$ fixed. For sufficiently
small $\tau>0$, Lemma~\ref{lem:bad} gives nonzero points and
steps, with $M_2\succ0$. Put
$s_1:=z_2-z_1$, $s_2:=z_3-z_2$ and
$\widehat\theta:=\lambda_Q(z_2,M_2)>0$.
The secant equation for $z_1\to z_2$ and stationarity at $z_2$ and $z_3$ give
\[
M_2s_1=s_1,\qquad z_2^Ts_1=0,\qquad
s_2=-\widehat\theta M_2^{-1}z_2,\qquad z_3^Ts_2=0.
\]
Equation~\eqref{eq:conjugacy}, applied to $z_1\to z_2\to z_3$,
gives $s_1^Ts_2=0$. Since $z_3=z_2+s_2$, both $z_2,z_3$ belong
to $s_1^\perp$.
This plane is $M_2$-invariant: for $u\perp s_1$,
$s_1^TM_2u=(M_2s_1)^Tu=s_1^Tu=0$.
Its normal eigenvalue is exactly $1$.

The resulting orthonormal frame is
\begin{equation}\label{eq:terminal-frame}
\left(-\frac{s_2}{\norm{s_2}},\,
       \frac{z_3}{\norm{z_3}},\,
       \frac{s_1}{\norm{s_1}}\right).
\end{equation}
In this frame, $z_2=(\norm{s_2},\norm{z_3},0)^T$ and
$z_3=(0,\norm{z_3},0)^T$. Since $Z_Q(z_2,M_2)=z_3\ne0$,
Lemma~\ref{lem:canonical} applies to this actual pair. It supplies
the following positive parameters; the stronger inequality for
$\widehat a$ is verified immediately below:
\begin{equation}\label{eq:terminal}
\begin{gathered}
\widehat L:=\norm{z_3-z_2}>0,\qquad
\widehat R:=\norm{z_3}>0,\qquad
\widehat\theta:=\lambda_Q(z_2,M_2)>0,\\
\widehat a:=\frac{z_3^TM_2z_3}{\widehat R^2}
 -\widehat\theta(\widehat R/\widehat L)^2
 >1+(\widehat R/\widehat L)^2.
\end{gathered}
\end{equation}
The canonical representation is
\[
M_2=\frac{\widehat\theta}{\widehat L^2}z_2z_2^T
     +\diag(0,\widehat a,1),
\]
whose determinant is $\det M_2=\widehat\theta\widehat a$.
Using \eqref{eq:z2}, \eqref{eq:z3} and \eqref{eq:M2-limit} therefore gives
\[
\frac{\widehat R}{\widehat L}\to\frac{b-1}{b+1},\qquad
\widehat\theta\to\frac{b(b+1)}{b^2+1},\qquad
\widehat a=\frac{\det M_2}{\widehat\theta}
 \to\frac{b^2+1}{b+1}.
\]
Consequently,
\[
\widehat a-1-(\widehat R/\widehat L)^2
\longrightarrow\frac{(b-1)(b^2+1)}{(b+1)^2}>0,
\]
which proves the strict margin in \eqref{eq:terminal} for all
sufficiently small $\tau>0$.

Keeping $g_3=z_3$ would give, by Lemma~\ref{lem:endpoint-update},
the auxiliary matrix
\[
M_3=\U(M_2;s_2,s_2)=\diag(1,\widehat a,1)
\]
in the terminal frame. Thus $z_3=\widehat R e_2$ is its eigenvector,
$\lambda_Q(z_3,M_3)=\widehat a$ and $Z_Q(z_3,M_3)=0$.
We must change this gradient to obtain a nonzero next point. Fix one such $\tau$ from now on; the points, $M_2$, the
terminal parameters and the current quotient are all fixed.

For the local calculation only, write $L,R,a$ for
$\widehat L,\widehat R,\widehat a$; thus $a>1+(R/L)^2$.
Work in \eqref{eq:terminal-frame}, and denote the coordinates of
the fixed point $z_3$ in this frame by $\widetilde z_3:=(0,R,0)^T$.
The parameterized output $z_{4,\epsilon}(t)$ below is also expressed
in this terminal frame.
For a size parameter $\epsilon$ and an initially free direction
parameter $t$, define the actual local data
\[
\begin{gathered}
g_{\epsilon,t}:=(0,R+\epsilon t,\epsilon)^T,\qquad
M_{\epsilon,t}:=\mathcal M_{L,a,1}(\epsilon t,\epsilon),\\
(z_{4,\epsilon}(t),D_\epsilon(t))
 :=\Sstep(\widetilde z_3,g_{\epsilon,t},M_{\epsilon,t}).
\end{gathered}
\]
Lemma~\ref{lem:endpoint-update} identifies $M_{\epsilon,t}$ as the
incoming update for the fixed displacement $-Le_1$, with curvature
$L^2$. Also
\[
\widetilde z_3^TM_{\epsilon,t}^{-1}g_{\epsilon,t}
 =R(R+\epsilon t)/a>0
\]
for small $|\epsilon|$ and bounded $t$. Lemma~\ref{lem:transition}
therefore gives a well-defined outgoing step, positive curvature
and $D_\epsilon(t)\succ0$.

Let $O_{\rm term}$ be
the orthogonal matrix whose columns are the axes in
\eqref{eq:terminal-frame}, expressed in current cycle coordinates.
It maps terminal coordinates back to current cycle coordinates.
For the free direction parameter $t$, transport this family by
\[
\begin{gathered}
g_3:=O_{\rm term}g_{\epsilon,t},\qquad
B_{\rm reg}:=O_{\rm term}M_{\epsilon,t}O_{\rm term}^T,\\
z_4:=O_{\rm term}z_{4,\epsilon}(t),\qquad
D_{\rm new}:=O_{\rm term}D_\epsilon(t)O_{\rm term}^T.
\end{gathered}
\]
Equivalently,
\[
g_3=z_3+\epsilon\left(
 t\frac{z_3}{\norm{z_3}}
 +\frac{z_2-z_1}{\norm{z_2-z_1}}\right).
\]
The local incoming update and Lemma~\ref{lem:transition}'s
orthogonal covariance give the actual identities
\[
\begin{gathered}
g_3^Ts_2=0,\qquad
s_2^T(g_3-z_2)=\widehat L^2>0,\\
B_{\rm reg}=\U(M_2;s_2,g_3-z_2)\succ0,\\
\Sstep(z_3,g_3,B_{\rm reg})=(z_4,D_{\rm new}).
\end{gathered}
\]
The fifth step therefore has positive line parameter and curvature.

We select $t$ and verify the three conditions in
Definition~\ref{def:cycle-input} for cycle $m+1$.
\begin{enumerate}
\item The actual call above prescribes gradient
$z_4$ and gives $D_{\rm new}=D_{\rm new}^T\succ0$ by
Lemma~\ref{lem:transition}. For nonvanishing, we use the local expansion.

Section~\ref{sec:bridge-calculations} proves analyticity of these outputs,
with $z_{4,0}(t)=0$, $D_0(t)=I$, and the expansions
\begin{equation}\label{eq:variation}
z_{4,\epsilon}(t)=\epsilon v(t)+O(\epsilon^2),\qquad
D_\epsilon(t)=I+\epsilon E(t)+O(\epsilon^2),
\end{equation}
uniformly on compact $t$-intervals, where
\begin{equation}\label{eq:vE}
v(t):=\begin{pmatrix}-\eta t\\0\\-a\end{pmatrix},\qquad
E(t):=\begin{pmatrix}
0&0&-1/L\\
0&t/R&1/R\\
-1/L&1/R&0
\end{pmatrix},\qquad \eta:=R/L.
\end{equation}
Since $v_3(t)=-a\ne0$, the expansion gives $z_{4,\epsilon}(t)\ne0$
for small positive $\epsilon$, uniformly on compact $t$-intervals.
Thus $z_4\ne0$ for the bounded choice of $t$ made next.

\item We impose invariance on the actual outputs by setting, for
$\epsilon\ne0$,
\[
\mathcal D(\epsilon,t):=\epsilon^{-6}
 \det[z_{4,\epsilon}(t),D_\epsilon(t)z_{4,\epsilon}(t),
                   D_\epsilon(t)^2z_{4,\epsilon}(t)].
\]
Subtracting columns changes the numerator to
\[
\det[z_{4,\epsilon}(t),(D_\epsilon(t)-I)z_{4,\epsilon}(t),
                      (D_\epsilon(t)-I)^2z_{4,\epsilon}(t)].
\]
Because $z_{4,0}(t)=0$ and $D_0(t)=I$ for every nearby $t$, these
columns have analytic factors $\epsilon$, $\epsilon^2$ and
$\epsilon^3$, respectively. Hence $\mathcal D$ extends analytically
to $\epsilon=0$. The determinant calculation in
Section~\ref{sec:bridge-calculations} gives
\begin{equation}\label{eq:determinant}
\begin{aligned}
\mathcal D(0,t)&=\det[v(t),E(t)v(t),E(t)^2v(t)]\\
&=\frac{t}{\eta^2L^3}
 \left[a^2(a-1-\eta^2)-\eta^2t^2(a-\eta^2)\right].
\end{aligned}
\end{equation}
Put
\[
q:=\sqrt{\frac{a-1-\eta^2}{a-\eta^2}},\qquad
 t_0:=\frac{aq}{\eta}.
\]
The terminal margin gives $0<q<1$ and $t_0>0$.
Moreover, $\mathcal D(0,t_0)=0$ and
$\partial_t\mathcal D(0,t_0)=-2t_0^2(a-\eta^2)/L^3<0$.
The analytic implicit function theorem therefore supplies a unique
local branch satisfying
\begin{equation}\label{eq:root}
\mathcal D(\epsilon,t(\epsilon))=0,\qquad
 t(\epsilon)=t_0+O(\epsilon)>0.
\end{equation}
From now on evaluate all local and transported families at
$t=t(\epsilon)$, suppress that argument, and write
\[
g_\epsilon:=g_{\epsilon,t(\epsilon)},\quad
M_\epsilon:=M_{\epsilon,t(\epsilon)},\qquad
v:=v(t_0),\quad E:=E(t_0).
\]
In particular, $z_{4,\epsilon}=\epsilon v+O(\epsilon^2)$ and
$D_\epsilon=I+\epsilon E+O(\epsilon^2)$.
The limiting vectors
\[
\epsilon^{-1}z_{4,\epsilon}\to v,\qquad
\epsilon^{-2}(D_\epsilon-I)z_{4,\epsilon}\to Ev
\]
are independent, since $v_2=0$ and $(Ev)_2=-a/R\ne0$.
Thus $z_{4,\epsilon},D_\epsilon z_{4,\epsilon}$ are independent for
small positive $\epsilon$. The exact equation \eqref{eq:root}
now proves that
\[
P_\epsilon:=\Span(z_{4,\epsilon},D_\epsilon z_{4,\epsilon})
\]
is a two-dimensional $D_\epsilon$-invariant plane containing
$z_{4,\epsilon}$, which is not an eigenvector.

For the normal eigenvalue, use the actual unit normal
\[
n_\epsilon:=
 \frac{z_{4,\epsilon}\times D_\epsilon z_{4,\epsilon}}
      {\norm{z_{4,\epsilon}\times D_\epsilon z_{4,\epsilon}}}.
\]
The denominator is nonzero by independence. To compute its limit,
put $n:=(-1,\eta(q^2-1),q)^T$.
Multiplication using \eqref{eq:vE} and
$(a-\eta^2)(q^2-1)=-1$ gives
\[
v\times Ev=\frac{a^2}{\eta L}n,\qquad En=(q/L)n.
\]
Since
$z_{4,\epsilon}\times D_\epsilon z_{4,\epsilon}
=\epsilon^3(v\times Ev+O(\epsilon))$,
one has $n_\epsilon=n/\norm{n}+O(\epsilon)$ for $\epsilon>0$.
Symmetry and exact invariance make this vector a normal
eigenvector of the actual matrix, with eigenvalue
\begin{equation}\label{eq:normal-eigenvalue}
b_\epsilon:=n_\epsilon^TD_\epsilon n_\epsilon
 =1+\epsilon q/L+O(\epsilon^2)>1,\qquad b_\epsilon\to1.
\end{equation}

The transported plane $P_{\rm new}:=O_{\rm term}P_\epsilon$ is
therefore two-dimensional, contains $z_4$ and is exactly
$D_{\rm new}$-invariant. Its unit normal $O_{\rm term}n_\epsilon$
has eigenvalue $b_{\rm new}:=b_\epsilon$.

\item Define $z_{5,\epsilon}:=Z_Q(z_{4,\epsilon},D_\epsilon)$ and
$\mu:=v^TEv/\norm{v}^2$.
The Neumann expansion and \eqref{eq:variation} give
\[
D_\epsilon^{-1}z_{4,\epsilon}
 =z_{4,\epsilon}-\epsilon^2Ev+O(\epsilon^3).
\]
Since $v_3=-a\ne0$, $\norm{z_{4,\epsilon}/\epsilon}$ is bounded away
from $0$. Substitution in the line-parameter formula yields
\[
\lambda_Q(z_{4,\epsilon},D_\epsilon)
=\frac{\norm{z_{4,\epsilon}/\epsilon}^2
       -\epsilon v^TEv+O(\epsilon^2)}
      {\norm{z_{4,\epsilon}/\epsilon}^2
       -2\epsilon v^TEv+O(\epsilon^2)}
=1+\epsilon\mu+O(\epsilon^2).
\]
It follows that
\begin{equation}\label{eq:regen-endpoint}
z_{5,\epsilon}
=z_{4,\epsilon}-\lambda_Q(z_{4,\epsilon},D_\epsilon)D_\epsilon^{-1}z_{4,\epsilon}
=\epsilon^2(Ev-\mu v)+O(\epsilon^3).
\end{equation}
The leading vector has second component $-a/R\ne0$, so
$z_{5,\epsilon}\ne0$. Orthogonal covariance gives the next point in the original coordinates
\[
z_5:=O_{\rm term}z_{5,\epsilon}
=Z_Q(z_4,D_{\rm new})\ne0.
\]

It remains to check the bounds in \eqref{eq:cycle-state}.
The preceding expansions give
\[
\norm{z_{4,\epsilon}}=\epsilon\norm{v}+O(\epsilon^2),\qquad
\norm{z_{5,\epsilon}}=\epsilon^2\norm{Ev-\mu v}+O(\epsilon^3).
\]
Together with $\rho(1+h)=h/2+O(h^2)$, \eqref{eq:normal-eigenvalue}
therefore gives
\[
\epsilon\frac{\rho(b_\epsilon)}{\norm{z_{5,\epsilon}}}
 \longrightarrow\frac{q}{2L\norm{Ev-\mu v}}>0.
\]
We have proved all three limits
\begin{equation}\label{eq:capacity}
\norm{z_{4,\epsilon}}\to0,\qquad
\frac{\norm{z_{5,\epsilon}}}{\norm{z_{4,\epsilon}}}\to0,\qquad
\frac{\rho(b_\epsilon)}{\norm{z_{5,\epsilon}}}\to+\infty.
\end{equation}
Equation~\eqref{eq:normal-eigenvalue} gives $b_{\rm new}>1$ and
$b_{\rm new}\to1$, hence also $b_{\rm new}<11/10$ for small
positive $\epsilon$. The last two limits in \eqref{eq:capacity}
and orthogonality give
\[
\norm{z_5}<\norm{z_4}/8,\qquad
\rho(b_{\rm new})/\norm{z_5}\ge4(m+1)
\]
for all sufficiently small $\epsilon>0$.
\end{enumerate}

Setting $\delta:=g_3-z_3$ gives
\begin{equation}\label{eq:actual-auxiliary}
B_{\rm reg}-M_3
=\frac{s_2\delta^T+\delta s_2^T+\delta\delta^T}{\widehat L^2}.
\end{equation}
The two updates have the same subtraction term and curvature
$\widehat L^2$, with gradient differences $s_2+\delta$ and $s_2$.
The earlier points, $M_2$, terminal parameters and current quotient
remain fixed; the fifth step uses $B_{\rm reg}$ rather than $M_3$.

For the prescribed tolerance $\eta_*>0$, the limits
$\norm{g_3-z_3}=\norm{g_\epsilon-\widetilde z_3}
=\epsilon\sqrt{t(\epsilon)^2+1}\to0$ and
$\norm{z_4}\to0$ also give
\begin{equation}\label{eq:bridge-smallness}
\norm{g_3-z_3}\le\eta_* ,\qquad
0<\norm{z_4}<\norm{z_3}/8.
\end{equation}
Here $m$, $\eta_*$ and the terminal data, including
$\norm{z_3}>0$, are fixed. Each requirement holds for all
sufficiently small positive $\epsilon$, so finitely many positive
thresholds can be imposed simultaneously. Fix one such $\epsilon$.

To form the next block's direct input, Lemma~\ref{lem:canonical}
applied to the fixed data $(z_4,z_4,D_{\rm new})$ supplies
positive $L_{\rm new},R_{\rm new},\theta_{\rm new},a_{\rm new}$ and
the orthonormal frame
\begin{equation}\label{eq:bridge-frame}
\begin{gathered}
L_{\rm new}:=\norm{z_4-z_5},\qquad
R_{\rm new}:=\norm{z_5},\qquad
\theta_{\rm new}:=\lambda_Q(z_4,D_{\rm new}),\\
O_{\rm new}:=\left[
\frac{z_4-z_5}{L_{\rm new}},\quad
\frac{z_5}{R_{\rm new}},\quad O_{\rm term}n_\epsilon\right].
\end{gathered}
\end{equation}
Its denominators are positive by that lemma, since $z_5\ne0$ and $z_4-z_5
=\theta_{\rm new}D_{\rm new}^{-1}z_4\ne0$.
The matrix $O_{\rm new}$ maps the next adapted coordinates to
current cycle coordinates, with the exact identification
\begin{equation}\label{eq:bridge-canonical}
\begin{gathered}
O_{\rm new}^Tz_4=(L_{\rm new},R_{\rm new},0)^T,\qquad
O_{\rm new}^Tz_5=(0,R_{\rm new},0)^T,\\
O_{\rm new}^TD_{\rm new}O_{\rm new}
\text{ has form \eqref{eq:canonical} with }\\
(L,R,\theta,a,b)=(L_{\rm new},R_{\rm new},\theta_{\rm new},a_{\rm new},b_{\rm new}).
\end{gathered}
\end{equation}
The starting data for the next cycle are now fixed, with
$L_{\rm new},R_{\rm new},a_{\rm new}>0$ and $b_{\rm new}>1$.
The actual input of Lemma~\ref{lem:bad} is formed at $z_5$. For the next parameter $\tau_{\rm new}>0$, define
\begin{equation}\label{eq:bridge-input}
\begin{gathered}
\kappa_{\rm new}:=-\frac{b_{\rm new}^2R_{\rm new}^3}
 {a_{\rm new}(b_{\rm new}-1)L_{\rm new}^2},\qquad
T_{\rm new}:=\tau_{\rm new}R_{\rm new},\qquad
H_{\rm new}:=\kappa_{\rm new}\tau_{\rm new}^3,\\
z_0^{\rm new}:=O_{\rm new}^Tz_5=(0,R_{\rm new},0)^T,\qquad
g_0^{\rm new}:=(0,R_{\rm new}+T_{\rm new},H_{\rm new})^T,\\
M^{\rm new}:=\mathcal M_{L_{\rm new},a_{\rm new},b_{\rm new}}
(T_{\rm new},H_{\rm new}).
\end{gathered}
\end{equation}
Every denominator is positive. In particular,
$b_{\rm new}-1$ is fixed before $\tau_{\rm new}\downarrow0$,
so $\kappa_{\rm new}$ is finite.
Prescribe gradient $O_{\rm new}g_0^{\rm new}$ at $z_5$.
By \eqref{eq:bridge-canonical} and Lemma~\ref{lem:endpoint-update},
the incoming step has parameter $\theta_{\rm new}>0$, endpoint
stationarity and curvature $L_{\rm new}^2>0$. Its actual update,
transported by Lemma~\ref{lem:transition}, is
\[
O_{\rm new}^T\U(D_{\rm new};z_5-z_4,
                  O_{\rm new}g_0^{\rm new}-z_4)O_{\rm new}=M^{\rm new}\succ0.
\]
The first outgoing call is admissible by Lemma~\ref{lem:endpoint-update},
since $R_{\rm new}+T_{\rm new}=R_{\rm new}(1+\tau_{\rm new})>0$.
Thus the actual point, selected gradient and incoming matrix,
not just the scalar parameters, match Lemma~\ref{lem:bad}.
Its remaining admissibility, nonzero outputs and positive
curvatures hold for sufficiently small $\tau_{\rm new}>0$.
Lemma~\ref{lem:transition} transports these calls to current cycle
coordinates, preserving the step relations and Euclidean norms.

The simultaneous quantitative choice of $\tau_{\rm new}$ is made
in Section~\ref{sec:trajectory} for these fixed starting data. All
thresholds and asymptotic constants may depend on the data already
fixed; none is assumed uniform over infinitely many cycles.
\end{proof}

\subsection{Expansion and determinant calculation}\label{sec:bridge-calculations}

We verify \eqref{eq:variation} and \eqref{eq:determinant}, used
in the preceding proof. The expansions require only $L,R,a>0$;
the inequality $a>1+(R/L)^2$ supplies the positive simple root.

Use $g_{\epsilon,t}$, $M_{\epsilon,t}$ and
$(z_{4,\epsilon}(t),D_\epsilon(t))$ from the proof of
Lemma~\ref{lem:bridge}, with $t$ in a compact interval $I_0$.
For this expansion, define
\[
w(\epsilon,t):=M_{\epsilon,t}^{-1}g_{\epsilon,t},
\qquad
\lambda(\epsilon,t):=
\frac{\widetilde z_3^Tw(\epsilon,t)}{\norm{w(\epsilon,t)}^2},
\]
and let
\[
s(\epsilon,t):=z_{4,\epsilon}(t)-\widetilde z_3,
\qquad
y(\epsilon,t):=z_{4,\epsilon}(t)-g_{\epsilon,t}.
\]
At $\epsilon=0$, the input matrix is $M_{0,t}=\diag(1,a,1)$ and
\[
\begin{gathered}
w(0,t)=(R/a)e_2,\qquad \lambda(0,t)=a,\qquad z_{4,0}(t)=0,\\
D_0(t)=I,\qquad s(0,t)=y(0,t)=-Re_2.
\end{gathered}
\]
The projection denominator and the two update denominators are,
respectively, $R^2/a^2$, $aR^2$, and $R^2$. Also
$\widetilde z_3^TM_{\epsilon,t}^{-1}g_{\epsilon,t}=R(R+\epsilon t)/a>0$ uniformly
on $I_0$ for sufficiently small $|\epsilon|$. Thus all rational
formulas are well defined and analytic on a neighborhood of
$\{0\}\times I_0$.

Primes in this step denote $\partial_\epsilon$ at $\epsilon=0$,
with $t$ fixed; in particular,
$z_4':=\partial_\epsilon z_{4,\epsilon}(t)|_{\epsilon=0}$. Formula
\eqref{eq:inverse} gives
\[
w'=\begin{pmatrix}Rt/(aL)\\t/a\\1\end{pmatrix},
\qquad \lambda'=-at/R,
\qquad z_4'=-\lambda'w(0,t)-aw'=v(t).
\]

To find the matrix derivative, first use the secant equation
$D_\epsilon(t)s(\epsilon,t)=y(\epsilon,t)$, with
$y=s-(0,\epsilon t,\epsilon)^T$. Differentiating it at $0$ yields
\[
-RD'e_2+s'=s'-(0,t,1)^T,
\qquad D'e_2=(0,t/R,1/R)^T.
\]

For the remaining $(e_1,e_3)$ principal block, consider the
two-sided compression to $\Span(e_1,e_3)$; this computation does
not require that subspace to be invariant. The projections of both
$M_{\epsilon,t}s$ and $y$ vanish at $0$. Hence the restrictions of
the two rank-one terms in \eqref{eq:update} to this subspace are
$O(\epsilon^2)$. The corresponding block of $D'$ is
therefore the block of $M_{\epsilon,t}'$, namely
$\bigl(\begin{smallmatrix}0&-1/L\\-1/L&0\end{smallmatrix}\bigr)$.
Symmetry supplies the remaining entries, giving
$D'=E(t)$. On a smaller compact neighborhood of
$\{0\}\times I_0$ inside the analytic domain, the second
$\epsilon$ derivatives are uniformly bounded. Taylor's theorem
therefore proves \eqref{eq:variation} with the stated uniform remainders.

Multiplication in \eqref{eq:vE} gives
\[
E(t)v(t)=\frac1L
\begin{pmatrix}a\\-a/\eta\\\eta t\end{pmatrix},
\qquad
E(t)^2v(t)=\frac1{L^2}
\begin{pmatrix}
-\eta t\\t(\eta^2-a)/\eta^2\\-a(\eta^2+1)/\eta^2
\end{pmatrix}.
\]
Factoring $L^{-1}$ and $L^{-2}$ from the second and third columns
and expanding gives exactly \eqref{eq:determinant}.
At $t_0=aq/\eta$, the bracket in that formula is $0$; its
derivative is $-2\eta^2t_0(a-\eta^2)$. Thus the determinant
derivative at $t_0$ is $-2t_0^2(a-\eta^2)/L^3$, as used in the
exact-root argument.

\section{Proof of smoothness in Lemma~\ref{lem:realization}}\label{app:interpolation}

This appendix proves the analytic assertions about the objective
\eqref{eq:bumps} constructed in Section~\ref{sec:smooth-objective}:
smoothness and flatness at the accumulation point, compact support
of $F-Q$, interpolation, and the global Hessian bounds.
The identification with the actual BFGS run is proved
in Section~\ref{sec:identification}.
Here $\overline{\mathbb B}(u,r)$ is a closed Euclidean ball, and
$\mathrm D^\ell h$ is the $\ell$th Fr\'echet derivative, with
$\mathrm D^0h=h$. Derivative tensors carry the induced multilinear
operator norm, and
$\norm{\mathrm D^\ell h}_\infty=\sup_x\norm{\mathrm D^\ell h(x)}$.
\begin{proof}[Analytic part of Lemma~\ref{lem:realization}]

Define $\vartheta:\R\to\R$ and $\chi:\R^3\to\R$ by
\[
\vartheta(t):=
\begin{cases}e^{-1/t},&t>0,\\0,&t\le0,\end{cases}
\qquad
\chi(z):=\frac{\vartheta(1-\norm{z}^2)}
{\vartheta(1-\norm{z}^2)+\vartheta(\norm{z}^2-1/4)}.
\]
Every right derivative of $\vartheta$ at $0$ is the limit of a
polynomial in $1/t$ times $e^{-1/t}$ and hence is $0$; thus
$\vartheta\in C^\infty(\R)$. The denominator defining $\chi$ is
positive because its first term is positive for $\norm{z}<1$ and its
second for $\norm{z}>1/2$. Therefore $\chi\in C^\infty(\R^3)$,
$\chi(z)=1$ when $\norm{z}\le1/2$, and $\chi(z)=0$ when
$\norm{z}\ge1$.

For every unit vector $u\in\R^3$, set
$f_u(z):=(u^Tz)\chi(z)$. For each integer $\ell\ge0$, compact support
and continuity in $(u,z)$ give the finite constant
\[
K_\ell:=\sup_{\substack{u,z\in\R^3\\\norm{u}=1}}
\norm{\mathrm D^\ell f_u(z)}<\infty.
\]
For $\ell=2$, this is the constant used to fix $c_0$ in
Section~\ref{sec:smooth-objective}. Use $r_j,h_j,\delta_j$ and
the corrections $\psi_j$ defined there in \eqref{eq:bumps}.
The support of $\psi_j$ is contained in
$\overline{\mathbb B}(x_j,h_j)$. If $i<j$, then $r_j<r_i/8$ and
\[
\norm{x_i-x_j}\ge r_i-r_j>(r_i+r_j)/4=h_i+h_j.
\]
Thus these closed supporting balls are pairwise disjoint. On the
$j$th ball,
$3r_j/4\le\norm{x}\le5r_j/4$. The balls accumulate only at $0$, are
locally finite in $\R^3\setminus\{0\}$, and their union has closure in
$\overline{\mathbb B}(0,5r_1/4)$. Hence the sum in \eqref{eq:bumps} is
pointwise well defined, vanishes at $0$, and is smooth away from $0$.

We control derivatives of every order relative to the distance from $0$.

If $\delta_j\ne0$, scaling $f_u$ with
$u:=\delta_j/\norm{\delta_j}$ gives
\begin{equation}\label{eq:bump-bound}
\norm{\mathrm D^\ell\psi_j}_\infty
\le K_\ell\norm{\delta_j}h_j^{1-\ell}
\le K_\ell4^{\ell-1}c_0\,2^{-j}r_j^{j+4-\ell}.
\end{equation}
The same bound is immediate when $\delta_j=0$. Fix integers
$\ell,N\ge0$. Since $0<r_j<1$, for
$j\ge\max\{1,\ell+N-4\}$ and $x$ in the $j$th supporting ball,
\[
\frac{\norm{\mathrm D^\ell\psi_j(x)}}{\norm{x}^N}
\le K_\ell4^{\ell-1}c_0(4/3)^N\,2^{-j}
 r_j^{j+4-\ell-N}
\le K_\ell4^{\ell-1}c_0(4/3)^N\,2^{-j}\longrightarrow0.
\]
The finitely many earlier balls are absent from a sufficiently small
neighborhood of $0$, and all derivatives vanish outside the supporting
balls. Therefore this decay holds as $x\to0$, not merely along the
interpolation points.

We turn the decay estimates into successive Fr\'echet derivatives at $0$, not merely estimates
along the points $x_j$.

For every integer $\ell\ge0$, define the tensor field
\[
A_\ell(x):=
\begin{cases}
\sum_{j=1}^\infty\mathrm D^\ell\psi_j(x),&x\ne0,\\
0,&x=0.
\end{cases}
\]
The preceding bounds make $A_\ell$ continuous at $0$. Taking $N=2$
also gives $\norm{A_\ell(x)}/\norm{x}\to0$, so $A_\ell$ is
Fr\'echet differentiable at $0$ with derivative
$0=A_{\ell+1}(0)$. Away from $0$, local finiteness gives
$\mathrm DA_\ell=A_{\ell+1}$. Induction on $\ell$ proves that $F-Q$
is $C^\infty$ and that all its derivatives vanish at $0$. Its support
is compact because it is contained in
$\overline{\mathbb B}(0,5r_1/4)$.

Disjoint supports reduce the Hessian estimate to one correction at a time.

At any point, at most one bump derivative contributes. By
\eqref{eq:bump-bound} with $\ell=2$, and because
$2^{-j}r_j^{j+2}<1$,
\[
\norm{\mathrm D^2(F-Q)}_\infty
\le4K_2c_0
=\frac{K_2}{2(1+K_2)}<\frac12.
\]
Since $\nabla^2Q=I$, this proves the Hessian bounds in
\eqref{eq:realization}. In a neighborhood of $x_j$, one has
$\chi((x-x_j)/h_j)=1$, while every other bump vanishes. Hence
$F(x_j)=Q(x_j)$ and
$\nabla F(x_j)=x_j+\delta_j=g_j$. The assertions at $0$ follow from
flatness.

\end{proof}
\end{document}